\documentclass[12pt]{article}
\usepackage{amssymb}
\usepackage{amsthm}
\usepackage{amsmath}
\usepackage{graphicx}
 \newtheorem{thm}{Theorem}[section]
 
 \newtheorem{lem}[thm]{Lemma}
 
 \theoremstyle{definition}
 
 \newtheorem{rem}[thm]{Remark}
 \numberwithin{equation}{section}

\begin{document}
\title{Cauchy problem of nonlinear Schr${\rm \ddot{o}}$dinger equation with initial data
in Sobolev space $W^{s,p}$ for $p<2$}
\author{Yi
Zhou\thanks{School of Mathematics, Fudan University, Shanghai
200433, P. R. China; Key Laboratory of Mathematics for Nonlinear
Sciences (Fudan University), Ministry of Education, P. R. China,
{\it email: yizhou@fudan.ac.cn.}}}
\date{\today}
\maketitle

\begin{abstract} In this paper, we consider in $R^n$ the Cauchy problem for
nonlinear Schr${\rm \ddot{o}}$dinger equation with initial data in
Sobolev space $W^{s,p}$ for $p<2$. It is well known that this
problem is ill posed. However, We show that after a  linear
transformation by the linear semigroup the problem becomes locally
well posed in $W^{s,p}$ for  $\frac{2n}{n+1}<p<2$ and
$s>n(1-\frac{1}{p})$. Moreover, we show that in one space
dimension, the problem is locally well posed in $L^p$ for any
$1<p<2$.
\end{abstract}
\textbf{Keyword:}Cauchy problem, nonlinear Schr${\rm
\ddot{o}}$dinger equation, locally well-posedness, scaling limit.

\section{Introduction}
Consider the Cauchy problem for the linear Schr${\rm
\ddot{o}}$dinger equation
\begin{equation}\label{a1}
iu_t(t,x)-\triangle u(t,x)=0,
\end{equation}
\begin{equation}\label{a2}
u(0,x)=u_0(x),
\end{equation}
 where $\triangle$ is the Laplace operator in
$R^n$ for $n\ge 1$. It is well known that this problem is well
posed for initial data $u_0\in L^p(R^n)$ if and only if $p=2$. For
this reason, it is believed that one can not solve the nonlinear
Schr${\rm \ddot{o}}$dinger equation with initial data in the
Sobolev space $W^{s,p}$ for $p\not=2$. However, this is not quite
right.

Notice that the solution of the Cauchy problem for
\eqref{a1},\eqref{a2} can be written as
\begin{equation}\label{a3}
u(t)=S(t)u_0=E(t)*u_0,
\end{equation}
where
\begin{equation}\label{a4}
E(t,x)=\frac{1}{(-4\pi it)^\frac{n}{2}}e^{-i\frac{|x|^2}{4t}}
\end{equation}
is the fundamental solution and $S(t)$ defines a semigroup. Thus
\begin{equation}
S(-t)u(t)\equiv u_0,
\end{equation}
and for any norm $X$ we have
\begin{equation}\label{a5}
\|S(-t)u(t)\|_X=\|u_0\|_X.
\end{equation}
There are some examples in the literature that one studies the
nonlinear Schr${\rm \ddot{o}}$dinger equation by using the norm
\begin{equation}\label{a6}
\|u\|_Y\triangleq\|S(-t)u(t)\|_X
\end{equation}
where $X$ is the usual Sobolev or weighted Sobolev norm. Of
course, we have the trivial example that when $X=H^s$, we have
$X=Y$. The first nontrivial example is to take $X$ to be the
weighted $L^2$ norm. Thus, we take
\begin{equation}\label{a7}
\|w\|_X=\sum_{|\alpha|\le s}\|x^\alpha w\|_{L^2(R^n)}
\end{equation}
where $\alpha$ is a multi-index.
Then
\begin{equation}\label{a8}
\|u(t)\|_Y=\sum_{|\alpha |\le s}\|x^\alpha S(-t)u(t)\|_{L^2(R^n)}
=\sum_{|\alpha |\le s}\|S(t)x^\alpha S(-t)u(t)\|_{L^2(R^n)}.
\end{equation}
Noting that
\begin{equation}\label{a9}
S(t)x_kS(-t)=x_k-2it\partial_{x_k}\triangleq L_k,
\end{equation}
we obtain
\begin{equation}\label{a10}
\|u(t)\|_Y=\sum_{|\alpha |\le s}\|L^\alpha u(t)\|_{L^2(R^n)}.
\end{equation}
 This norm was
first used by McKean and Shatah \cite{McKeanShatah} and it was
proved that one has the following global Sobolev inequality
\begin{equation}\label{a13}
\|u(t)\|_{L^\infty}\le C(1+t)^{-\frac{n}{2}}(\sum_{|\alpha|\le
s}\|L^{\alpha }u(t)\|_{L^2(R^n)}+\|u(t)\|_{H^s(R^n)}),\quad
s>\frac{n}{2}.
\end{equation}
This inequalty is similar to the global Sobolev inequality for the
wave equation obtained earlier by Klainerman (see
\cite{Klainerman1}) and is very important in studying the
nonlinear problem in their paper.

Another more recent example is to take $X=H^b_tH^s_x$, then $Y$ is
the so called Bourgain space (see \cite{Bourgain}). This space
plays a very important role in the recent study of low regularity
solution of nonlinear Schr${\rm \ddot{o}}$dinger equations.

Therefore, why not take $X=L^p$ (or $W^{s,p}$)? It is our aim to
investigate this problem in this paper.

Consider the Cauchy problem for the nonlinear Schr${\rm
\ddot{o}}$dinger equation
\begin{equation}\label{a14}
iu_t(t,x)-\triangle u(t,x)=\pm |u(t,x)|^2u(t,x),
\end{equation}
\begin{equation}\label{a15}
u(0,x)=u_0(x).
\end{equation}
This problem can be reformulated as
\begin{equation}\label{a16}
u(t)=S(t)u_0\pm \int_0^tS(t-\tau )(|u(\tau )|^2u(\tau ))d\tau .
\end{equation}
Motivated by our above discussions, we make a linear
transformation
\begin{equation}\label{a17}
v(t)=S(-t)u(t),
\end{equation}
then
\begin{equation}\label{a18}
u(t)=S(t)v(t).
\end{equation}
Therefore, we get
\begin{equation}\label{a19}
v(t)=u_0\pm \int_0^tS(-\tau)[S(-\tau)\bar
v(\tau)(S(\tau)v(\tau))^2]d\tau ,
\end{equation}
where we use the fact that $\bar{S}(\tau)=S(-\tau )$.

Our main result in this paper is that \eqref{a19} is locally well
posed in Sobolev space $W^{s,p}$ for certain $p<2$. More
precisely, we have the following:
\begin{thm}Consider the nonlinear integral equation \eqref{a19},
suppose that
\begin{equation}\label{a20}
u_0\in W^{s,p}(R^n)
\end{equation}
for $s>n(1-\frac{1}{p})$ and $\frac{2n}{n+1}<p<2$, where
$W^{s,p}(R^n)$ is understood as $B^s_{p,p}(R^n)$ and
$B^s_{p,q}(R^n)$ is the Besov space. Then there exists a time T
which only depends on $\|u_0\|_{W^{s,p}(R^n)}$ such that the
integral equation has a unique solution $v\in C([0,T],
W^{s,p}(R^n))$  satisfying
\begin{equation}\label{a21}
\|v(t)\|_{W^{s,p}(R^n)}\le 2\|u_0\|_{W^{s,p}(R^n)},\quad \forall
t\in [0,T].
\end{equation}
 Moreover, suppose that $v_1$, $v_2$ are two solutions with
initial data $u_{01}$, $u_{02}$, then there holds
\begin{equation}\label{a22}
\|v_1(t)-v_2(t)\|_{W^{s,p}(R^n)}\le
2\|u_{01}-u_{02}\|_{W^{s,p}(R^n)},\quad \forall t\in [0,T].
\end{equation}
\end{thm}

\begin{rem}Our proof relays on a subtle cancellation in the
nonlinearity and thus our result is not valid for the general
nonlinearity $F(u,\bar u)$. However, for nonlinear term of the
form $\pm |u|^{2m}u$, where $m$ is an integer, it is not
difficulty to generalize our result to this case.
\end{rem}

\begin{rem}Similar results are expected for other nonlinear  dispersive
equations and nonlinear wave equations. However, no such result is
presently known.
\end{rem}

 We point out that we can also
slightly improve our result by using Besov spaces.
\begin{thm}Consider the nonlinear integral equation \eqref{a19},
suppose that
\begin{equation}\label{a26}
u_0\in \dot{B}^s_{p,1}(R^n)
\end{equation}
for $s=n(1-\frac{1}{p})$ and $\frac{2n}{n+1}<p<2$, where
$\dot{B}^s_{p,1}(R^n)$ is the homogenous Besov space. Then there
exists a time T which only depends on
$\|u_0\|_{\dot{B}^s_{p,1}(R^n)}$ such that the integral equation
has a unique solution $v\in C([0,T], \dot{B}^s_{p,1}(R^n))$
satisfying
\begin{equation}\label{a27}
\|v(t)\|_{\dot{B}^s_{p,1}(R^n)}\le
2\|u_0\|_{\dot{B}^s_{p,1}(R^n)}, \quad \forall t\in [0,T].
\end{equation}
 Moreover, suppose that $v_1$, $v_2$ are two solutions with
initial data $u_{01}$, $u_{02}$, then there holds
\begin{equation}\label{a28}
\|v_1(t)-v_2(t)\|_{\dot{B}^s_{p,1}(R^n)}\le
2\|u_{01}-u_{02}\|_{\dot{B}^s_{p,1}(R^n)},\quad \forall t\in
[0,T].
\end{equation}

\end{thm}
In the following, we will only prove Theorem 1.4 since the proof
of Theorem 1.1 is similar.

We point out that Theorem 1.1 is only to show that one can solve
the Cauchy problem in $W^{s,p}$ for $p<2$, the regularity
assumption in Theorem 1.1 need not be optimal and can be improved.
As an example, we will show that the problem is locally well posed
in $L^p$ for any $1<p<2$ in one space dimension. It is proved by
Y.Tsutsumi \cite{Tsutsumi} that the problem is locally well-posed
in $L^2$.Then It is proved by Gr${\rm \ddot{u}}$nrock
\cite{Gruenrock} that the problem is locally well posed in
$\hat{L}^p$, for any $1<p<\infty$ (se also Cazenave et al
\cite{Cazenave et al} and Vargas and Vega \cite{VargasVega}.) Here
\begin{equation}
\|f\|_{\hat{L}^p}=\|\hat{f}\|_{L^{p'}},
\end{equation}
where $\hat{f}$ is the Fourier transform of f and
\begin{equation}
\frac{1}{p}+\frac{1}{p'}=1.
\end{equation}
Noting that
\begin{equation}
\|\hat{f}\|_{L^{p'}}\le C\|f\|_{L^p},\quad 1\le p\le 2,
\end{equation}
$\hat{L}^p$ space is slightly larger than $L^p$ space.
However,$L^p$ is more commonly used space. More recently, there
are even some local existence result in $H^s$ for some $s<0$, see
Christ et al \cite{Christ et al} as well as Koch and Tataru
\cite{KochTataru}.

Our main result in one space dimension is as follows:
\begin{thm}Consider the nonlinear integral equation \eqref{a19} in
one space dimension, suppose that
\begin{equation}u_0\in L^p(R)
\end{equation}
for $1<p<2$. Then there exists a time T which only depends on
$\|u_0\|_{L^p(R)}$ such that the integral equation has a unique
solution $v\in C([0,T],L^p(R))$  satisfying
\begin{equation}\label{a29}
\|v(t)\|_{L^p(R)}\le C_0\|u_0\|_{L^p(R)},\quad\forall t\in [0,T],
\end{equation}
and
\begin{equation}\label{a80}
\left\{\int_0^T\tau^{\theta p'}\|\partial_\tau v(\tau
)\|_{L^p(R)}^{p'}d\tau \right\}^{\frac{1}{p'}}\le
C_1\|u_0\|_{L^p(R)}^3.
\end{equation}
where
\begin{equation}\label{d11}
\frac{1}{p}+\frac{1}{p'}=1,\quad \theta =\frac{2}{p}-1.
\end{equation}
 Moreover, suppose that $v_1$, $v_2$ are two solutions with
initial data $u_{01}$, $u_{02}$, then there holds
\begin{equation}\label{a30}
\|v_1(t)-v_2(t)\|_{L^p(R)}\le
C_0\|u_{01}-u_{02}\|_{L^p(R)},\quad\forall t\in [0,T].
\end{equation}
Here $C_0$ and $C_1$ are positive constants independent of the
initial data.
\end{thm}
\begin{rem} Let $u(t,x)$ be a solution to the nonlinear
Schr${\rm \ddot{o}}$dinger equation \eqref{a14} with initial data
\eqref{a15}, then $u_{\lambda}(t,x)=\lambda u(\lambda^2 t,\lambda
x)$ is also a solution with initial data $u_{0\lambda}=\lambda
u_0(\lambda x)$. If
\begin{equation}\|u_{0\lambda }\|_{L^p(R^n)}\equiv
\|u_0\|_{L^p(R^n)},
\end{equation}
then $p$ is called a scaling limit. It is easy to see that $p$ is
a scaling limit in one space dimension if and only if $p=1$. Thus,
as $p$ close to 1, we can go arbitrary close to the scaling limit.
\end{rem}

In the following, $C$ will denote a positive constant independent
of the initial data and its meaning may change from line to line.

 Finally, we refer to \cite{Bergh} for the definition of Besov spaces and homogenous Besov spaces.

\section{A Key Lemma}
A key Lemma leading to our local well posedness is the following:
\begin{lem}
We consider a trilinear form
\begin{equation}\label{b1}
v_0(\tau)=T(v_1(\tau), v_2(\tau ),
v_3(\tau))=S(-\tau)[S(-\tau)v_1(\tau)S(\tau )v_2(\tau )S(\tau
)v_3(\tau )].
\end{equation}
Then, there holds
\begin{equation}\label{b2}
\|v_0(\tau )\|_{L^1(R^n)}\le C\tau ^{-n}\|v_1(\tau
)\|_{L^1(R^n)}\|v_2(\tau )\|_{L^1(R^n)}\|v_3(\tau )\|_{L^1(R^n)}.
\end{equation}
\end{lem}
\begin{proof}We directly compute
\begin{equation}\label{b3}
v_0(\tau, \alpha )=C\tau^{-2n}\iiiint e^{i\frac{|\alpha -\beta
|^2+|\beta -x|^2-|\beta -y|^2-|\beta -z|^2}{4\tau }} v_1(\tau
,x)v_2(\tau ,y)v_3(\tau ,z) d\beta dxdydz.
\end{equation}
It is easy to see that
\begin{eqnarray}\label{b4}
&& \tau^{-n}\int e^{i\frac{|\alpha -\beta |^2+|\beta -x|^2-|\beta
-y|^2-|\beta -z|^2}{4\tau }}d\beta \\\nonumber &&=\tau^{-n}
e^{i\frac{|\alpha |^2+|x|^2-|y|^2-|z|^2}{4\tau }}\int
e^{-i\frac{(\alpha +x-y-z)\cdot \beta }{2\tau }} d\beta
\\\nonumber && =2^n e^{i\frac{|\alpha |^2+|x|^2-|y|^2-|z|^2}{4\tau
}}\int e^{-i(\alpha +x-y-z)\cdot \gamma  } d\gamma\\\nonumber&& =C
e^{i\frac{|\alpha |^2+|x|^2-|y|^2-|z|^2}{4\tau }} \delta (\alpha
+x-y-z)
\end{eqnarray}
where $\delta $ is the Dirac function. The rest part of the proof
is obvious.
\end{proof}

\begin{lem}Let $v_l, l=0,1,2,3$ satisfy \eqref{b1}, suppose that
$2^{j-2}\le |\xi|\le 2^{j+2}$ in the support of $\hat{v}_2(\tau,
\xi)$ and $2^{k-2}\le |\xi|\le 2^{k+2}$ in the support of
$\hat{v}_3(\tau,\xi)$, where $\hat{v}_2$, $\hat{v}_3$ denote the
space Fourier transform of $v_2$,$v_3$. Then there holds
\begin{equation}\label{b7}
\|v_0(\tau )\|_{L^2}\le C2^{\frac{n}{2}(j+k)}\|v_1(\tau
)\|_{L^2}\|v_2(\tau )\|_{L^2}\|v_3(\tau )\|_{L^2}.
\end{equation}
\end{lem}
\begin{proof}Let $u_l(\tau )=S(\tau )v_l(\tau), l=0,2,3$ and $u_1(\tau) =S(-\tau )v_1(\tau
)$, then $ \hat{u}_0(\tau ,\xi)=e^{-i|\xi |^2\tau }\hat{v}_0(\tau
,\xi)$ etc. We have
\begin{equation}\label{b9}
u_0(\tau )=u_1(\tau )u_2(\tau)u_3(\tau).
\end{equation}
Therefore
\begin{eqnarray}\label{b10}
&&\|v_0(\tau )\|_{L^2(R^n)}=\|u_0(\tau )\|_{L^2(R^n)}\le
\|u_1(\tau )\|_{L^2(R^n)}\|u_2(\tau )\|_{L^\infty (R^n)}\|u_3(\tau
)\|_{L^\infty (R^n)}\\\nonumber && \le C\|u_1(\tau
)\|_{L^2(R^n)}\|\hat{u}_2(\tau )\|_{L^1 (R^n)}\|\hat{u}_3(\tau
)\|_{L^1 (R^n)}\\\nonumber &&=C\|v_1(\tau
)\|_{L^2(R^n)}\|\hat{v}_2(\tau )\|_{L^1 (R^n)}\|\hat{v}_3(\tau
)\|_{L^1 (R^n)}.
\end{eqnarray}
Noting the support property of $\hat{v}_2(\tau )$ and
$\hat{v}_3(\tau )$,  the desired conclusion follows from Schwartz
inequality.
\end{proof}

By the interpolation theorem on the multi-linear functionals (see
\cite{Bergh} page 96 Theorem 4.4.1), we can interpolate the
inequality in Lemma 2.1 and Lemma 2.2 to get the following:
\begin{lem}Let $v_l, l=0,1,2,3$ satisfy \eqref{b1}. Suppose that
$2^{j-2}\le |\xi|\le 2^{j+2}$ in the support of $\hat{v}_2(\tau,
\xi)$ and $2^{k-2}\le |\xi|\le 2^{k+2}$ in the support of
$\hat{v}_3(\tau,\xi)$, where $\hat{v}_2$, $\hat{v}_3$ denote the
space Fourier transform of $v_2$,$v_3$. Then there holds
\begin{equation}\label{b12}
\|v_0(\tau )\|_{L^p}\le
C\tau^{-n(\frac{2}{p}-1)}2^{n(1-\frac{1}{p})(j+k)}\|v_1(\tau
)\|_{L^p}\|v_2(\tau )\|_{L^p}\|v_3(\tau )\|_{L^p},\quad 1\le p\le
2.
\end{equation}
\end{lem}

\section{Proof of the Theorem 1.4 }
In this section, we prove Theorem 1.4 by a contraction mapping
principle.

Let us define the space
\begin{equation}\label{c1} X=\{w\in C([0,T],
\dot{B}^s_{p,1}(R^n))| \sup_{0\le t\le
T}\|w(t)\|_{\dot{B}^s_{p,1}(R^n)}\le
2\|u_0\|_{\dot{B}^s_{p,1}(R^n)}\},
\end{equation}
where $s=n(1-\frac{1}{p})$ and $\frac{2n}{n+1}<p<2$. For any $w\in
X$, define a map M by
\begin{equation}\label{c2}
(Mw)(t)\triangleq u_0\pm \int_0^tS(-\tau)[S(-\tau)\bar
w(\tau)(S(\tau)w(\tau))^2]d\tau.
\end{equation}
We want to show that M maps X into itself and is a contraction
provided that T is sufficiently small.

Firstly let us recall the definition of homogenous Besov spaces.
Let $\psi \in C_0^{\infty}(R^n)$ such that
\begin{equation}\label{c3}
supp \psi \subset \{ \xi | |\xi |\le 1\}
\end{equation}
and
\begin{equation}\label{c4}
\psi (\xi )\equiv 1\quad |\xi |\le\frac{1}{2}.
\end{equation}
Let
\begin{equation}\label{c5}
\phi (\xi)=\psi (2^{-1}\xi)-\psi (\xi )
\end{equation}
then
\begin{equation}\label{c6}
\sum_{j=-\infty }^{+\infty }\phi (2^{-j}\xi )\equiv 1,
\end{equation}
 and we have the following dyadic decomposition
 \begin{equation}\label{c7}
 w=\sum_{j=-\infty}^{+\infty}w_j,
 \end{equation}
 where
 \begin{equation}\label{c8}
 \hat{w}_j(\xi )=\phi (2^{-j}\xi )\hat{w}(\xi ).
 \end{equation}
 The Besov norm $\dot {B}^s_{p,1}(R^n)$ is defined by
 \begin{equation}\label{c9}
 \|w\|_{\dot
 {B}^s_{p,1}(R^n)}=\sum_{j=-\infty}^{+\infty}2^{js}\|w_j\|_{L^p(R^n)}.
 \end{equation}
 Let $w\in X$, to show M maps X into itself, we need to estimate
 the nonlinear term
 \begin{eqnarray}\label{c10}
&&F(\tau )=S(-\tau )[S(-\tau )\bar{w}(\tau )(S(\tau )w(\tau
))^2]\\\nonumber && =\sum_{j,k,l}S(-\tau )[S(-\tau )\bar{w}_j(\tau
)S(\tau )w_k(\tau ) S(\tau )w_l(\tau )].
\end{eqnarray}
To estimate F, we only need to estimate
\begin{equation}\label{c11}
F_1(\tau )=\sum_{j\ge k\ge l}S(-\tau )[S(-\tau )\bar{w}_j(\tau
)S(\tau )w_k (\tau )S(\tau )w_l(\tau )],
\end{equation}
all the other terms in the summation can be estimated in a similar
way.

By Lemma 2.3, we have
\begin{eqnarray}\label{c12}
&& \|F_1(\tau )\|_{\dot{B}^s_{p,1}(R^n)}\\\nonumber &&\le
\sum_{j=-\infty}^{+\infty}\|\sum_{k,l=-\infty }^jS(-\tau )[S(-\tau
)\bar{w}_j(\tau )S(\tau )w_k(\tau ) S(\tau )w_l(\tau )]\|_{\dot
{B}^s_{p,1}(R^n)}\\\nonumber && \le
\sum_{j=-\infty}^{+\infty}\sum_{m=-\infty}^{j+4}2^{ms}\|\phi
(2^{-m}D)\{\sum_{k,l=-\infty }^jS(-\tau )[S(-\tau )\bar{w}_j(\tau
)S(\tau )w_k (\tau )S(\tau )w_l(\tau )]\}\|_{L^p(R^n)}\\\nonumber
&& \le
\sum_{j=-\infty}^{+\infty}\sum_{m=-\infty}^{j+4}2^{ms}\|\sum_{k,l=-\infty
}^jS(-\tau )[S(-\tau )\bar{w}_j(\tau )S(\tau )w_k (\tau )S(\tau
)w_l(\tau )]\|_{L^p(R^n)}\\\nonumber && \le
C\sum_{j=-\infty}^{+\infty}2^{js}\|\sum_{k,l=-\infty }^jS(-\tau
)[S(-\tau )\bar{w}_j(\tau )S(\tau )w_k (\tau )S(\tau )w_l(\tau
)]\|_{L^p(R^n)}\\\nonumber && \le C\sum_{j,k,l}2^{js}\|S(-\tau
)[S(-\tau )\bar{w}_j(\tau )S(\tau )w_k (\tau )S(\tau )w_l(\tau
)]\|_{L^p(R^n)}\\\nonumber && \le
C\tau^{-n(\frac{2}{p}-1)}\sum_{j,k,l}2^{(j+k+l)s}\|w_j(\tau
)\|_{L^p(R^n)}\|w_k(\tau )\|_{L^p(R^n)}\|w_l(\tau
)\|_{L^p(R^n)}\\\nonumber && =C\tau^{-n(\frac{2}{p}-1)}\|w(\tau
)\|_{\dot{B}^s_{p,1}(R^n)}^3,
\end{eqnarray}
where $s=n(1-\frac{1}{p})$. Therefore
\begin{equation}\label{c13}
 \|F(\tau )\|_{\dot{B}^s_{p,1}(R^n)}\le C\tau^{-n(\frac{2}{p}-1)}\|w(\tau
)\|_{\dot{B}^s_{p,1}(R^n)}^3.
\end{equation}
Noting that when $\frac{2n}{n+1}<p<2$, we have
$0<n(\frac{2}{p}-1)<1$, it is easy to see
\begin{eqnarray}\label{c14}
&& \|(Mw)(\tau )\|_{\dot{B}^s_{p,1}(R^n)}\le
\|u_0\|_{\dot{B}^s_{p,1}(R^n)}+\int_0^t\|F(\tau)\|_{\dot{B}^s_{p,1}(R^n)}\\\nonumber
&& \le  \|u_0\|_{\dot{B}^s_{p,1}(R^n)}+
 C\int_0^t\tau^{-n(\frac{2}{p}-1)}\|w(\tau)\|_{\dot{B}^s_{p,1}(R^n)}^3d\tau
\\\nonumber &&\le  \|u_0\|_{\dot{B}^s_{p,1}(R^n)}+
 CT^{1-n(\frac{2}{p}-1)}(\sup_{0\le t\le
 T}\|w(t)\|_{\dot{B}^s_{p,1}(R^n)})^3\\\nonumber &&\le  \|u_0\|_{\dot{B}^s_{p,1}(R^n)}+
 CT^{1-n(\frac{2}{p}-1)}\|u_0\|_{\dot{B}^s_{p,1}(R^n)}^3\\\nonumber &&
 \le  2\|u_0\|_{\dot{B}^s_{p,1}(R^n)}
 \end{eqnarray}
 provided that T is sufficiently small.

 Now we prove that M is a contraction.Let $w^{(1)},w^{(2)}\in X$,
 denote $w^\ast =w^{(1)}-w^{(2)}$ and $v^\ast =Mw^{(1)}-Mw^{(2)}$,
 then
 \begin{eqnarray}\label{c15}
&& v^\ast=\pm\int_0^tS(-\tau)[S(-\tau)\bar{w}^{(1)}(\tau )(S(\tau
)w^{(1)}(\tau ))^2-S(-\tau)\bar{w}^{(2)}(\tau )(S(\tau
)w^{(2)}(\tau ))^2]d\tau\\\nonumber &&
=\pm\int_0^tS(-\tau)[S(-\tau)\bar{w}^\ast (\tau )(S(\tau
)w^{(1)}(\tau ))^2-S(-\tau)\bar{w}^{(2)}(\tau )S(\tau
)(w^{(1)}(\tau )+w^{(2)}(\tau ))S(\tau )w^\ast (\tau )]d\tau
\end{eqnarray}
By a similar argument as before, we can get
\begin{eqnarray}\label{c16}
&&\|v^\ast (t)\|_{\dot{B}^s_{p,1}(R^n)}\\\nonumber &&\le
C\int_0^t\tau^{-n(\frac{2}{p}-1)}
(\|w^{(1)}(\tau)\|_{\dot{B}^s_{p,1}(R^n)}+\|w^{(2)}(\tau)\|_{\dot{B}^s_{p,1}(R^n)})^2
\|w^\ast(\tau)\|_{\dot{B}^s_{p,1}(R^n)} \\\nonumber &&\le
CT^{1-n(\frac{2}{p}-1)}\|u_0\|_{\dot{B}^s_{p,1}(R^n)}^2\sup_{0\le
t\le T}\|w^\ast(t)\|_{\dot{B}^s_{p,1}(R^n)}\\\nonumber && \le
\frac{1}{2}\sup_{0\le t\le T}\|w^\ast(t)\|_{\dot{B}^s_{p,1}(R^n)}.
\end{eqnarray}
Therefore, we proved the existence and uniqueness of the solution.
To prove the stability result, let $v^{(1)}$ and $v^{(2)}$ be two
solutions with initial data $u_{01}$ and $u_{02}$. With a little
abuse of notation, we still denote $v^\ast =v^{(1)}-v^{(2)}$. Then
we have
 \begin{eqnarray}\label{c17}
&& v^\ast=u_{01}-u_{02}\\\nonumber
&&\pm\int_0^tS(-\tau)[S(-\tau)\bar{v}^{(1)}(\tau )(S(\tau
)v^{(1)}(\tau ))^2-S(-\tau)\bar{v}^{(2)}(\tau )(S(\tau
)v^{(2)}(\tau ))^2]d\tau\\\nonumber && =u_{01}-u_{02}\\\nonumber
&&\pm\int_0^tS(-\tau)[S(-\tau)\bar{v}^\ast (\tau )(S(\tau
)v^{(1)}(\tau ))^2-S(-\tau)\bar{v}^{(2)}(\tau )S(\tau
)(v^{(1)}(\tau )+v^{(2)}(\tau ))S(\tau )v^\ast (\tau )]d\tau.
\end{eqnarray}
Thus,
\begin{eqnarray}\label{c18}
&&\|v^\ast (t)\|_{\dot{B}^s_{p,1}(R^n)}\le
\|u_{01}-u_{02}\|_{\dot{B}^s_{p,1}(R^n)}\\\nonumber
&&+C\int_0^t\tau^{-n(\frac{2}{p}-1)}
(\|v^{(1)}(\tau)\|_{\dot{B}^s_{p,1}(R^n)}+\|v^{(2)}(\tau)\|_{\dot{B}^s_{p,1}(R^n)})^2
\|v^\ast(\tau)\|_{\dot{B}^s_{p,1}(R^n)} \\\nonumber &&\le
\|u_{01}-u_{02}\|_{\dot{B}^s_{p,1}(R^n)}+
CT^{1-n(\frac{2}{p}-1)}(\|u_{01}\|_{\dot{B}^s_{p,1}(R^n)}+\|u_{02}\|_{\dot{B}^s_{p,1}(R^n)})^2\sup_{0\le
t\le T}\|v^\ast(t)\|_{\dot{B}^s_{p,1}(R^n)}\\\nonumber && \le
\|u_{01}-u_{02}\|_{\dot{B}^s_{p,1}(R^n)}+ \frac{1}{2}\sup_{0\le
t\le T}\|v^\ast(t)\|_{\dot{B}^s_{p,1}(R^n)}.
\end{eqnarray}
Therefore
\begin{equation}\label{c19} \sup_{0\le t\le
T}\|v^\ast(t)\|_{\dot{B}^s_{p,1}(R^n)}\le
2\|u_{01}-u_{02}\|_{\dot{B}^s_{p,1}(R^n)}.
\end{equation}
We completed the proof of Theorem 1.4.

\section{Proof of the Theorem 1.5 }
In this section, we will prove theorem 1.5.
\begin{lem}Let $n=1$ and $v_l$, $l=0,1,2,3$ be defined by Lemma
2.1, then there holds
\begin{equation}\label{d1}
\sup_{0\le \tau\le T}(\tau\|v_0(\tau)\|_{L^1(R)})\le
C\prod_{i=1}^3\{\|v_i(0)\|_{L^1(R)}+\int_0^T\|\partial_\tau
v_i(\tau)\|_{L^1(R)}d\tau\}.
\end{equation}
\end{lem}
\begin{proof}  \eqref{d1} follows from Lemma 2.1 by
\begin{equation}
v_i(t )=v_i(0)+\int_0^t\partial_\tau v_i(\tau )d\tau .
\end{equation}
\end{proof}

\begin{lem}Let $n=1$ and $v_l$, $l=0,1,2,3$ be  defined by Lemma
2.1, then there holds
\begin{equation}\label{d3}
\left\{\int_0^T\|v_0(\tau
)\|_{L^2(R)}^2d\tau\right\}^{\frac{1}{2}}\le
C\prod_{i=1}^3\{\|v_i(0)\|_{L^2(R)}+\int_0^T\|\partial_\tau
v_i(\tau)\|_{L^2(R)}d\tau\}.
\end{equation}
\end{lem}
\begin{proof}Let
\begin{equation}
u_1(\tau)=S(\tau )\bar v_1(\tau ),\quad
u_2(\tau)=S(\tau)v_2(\tau),\quad u_3(\tau)=S(\tau)v_3(\tau),
\end{equation}
then it follows from H${\rm \ddot{o}}$lder's inequality that
\begin{equation}\left\{\int_0^T\|v_0(\tau
)\|_{L^2(R)}^2d\tau\right\}^{\frac{1}{2}}\le
C\prod_{i=1}^3\left\{\int_0^T\|u_i(\tau
)\|_{L^6(R)}^6d\tau\right\}^{\frac{1}{6}}.
\end{equation}
Noting that
\begin{equation}iu_{1t}(t,x)-\triangle u_1(t,x)=S(t)\partial_t\bar
v_1(t),
\end{equation}
\begin{equation}
u_1(0)=\bar{v}_1(0)
\end{equation}
as well as similar equations for $u_2$, $u_3$, the desired
conclusion follows from Strichartz' inequality.
\end{proof}

By the interpolation theorem on the multi-linear functionals (see
[1] page 96 Theorem 4.4.1), we can interpolate the inequality in
Lemma 4.1 and Lemma 4.2 to get the following
\begin{lem} Let n=1 and $v_l$ $l=0,1,2,3$ be defined by Lemma 2.1,
then there holds
\begin{equation}\label{d10}
\left\{\int_0^T\tau^{\theta p'}\|v_0(\tau )\|_{L^p(R)}^{p'}d\tau
\right\}^{\frac{1}{p'}}\le
C\prod_{i=1}^3\{\|v_i(0)\|_{L^p(R)}+\int_0^T\|\partial_\tau
v_i(\tau)\|_{L^p(R)}d\tau\},
\end{equation}
where $1<p<2$ and $p'$,$\theta$ satisfy \eqref{d11}.
\end{lem}

We are now ready to prove Theorem 1.5.

Let us define the set
\begin{equation}\label{d12}
X=\{w| w(0)=u_0,\left\{\int_0^T\tau^{\theta p'}\|\partial_\tau
w(\tau )\|_{L^p(R)}^{p'}d\tau \right\}^{\frac{1}{p'}}\le
C_1\|u_0\|_{L^p(R)}^3\}
\end{equation}
where $\theta$,$p'$ are defined by \eqref{d11} and $C_1$ is a
positive constant independent of the initial data and will be
determined later. For any $w\in X$, define a map M by
\begin{equation}
(Mw)(t)=u_0\pm \int_0^tS(-\tau )[S(-\tau )\bar{w}(\tau
)(S(\tau)w(\tau ))^2]d\tau .
\end{equation}
We want to show that M maps $X$ into itself and is a contraction.

For simplicity, we denote $v=Mw$. Obviously,
\begin{equation}
v(0)=u_0
\end{equation}
and
\begin{equation}
\partial_\tau v(\tau )=\pm S(-\tau )[S(-\tau )\bar{w}(\tau
)(S(\tau)w(\tau ))^2.
\end{equation}
Applying Lemma 4.3, we get,
\begin{equation}
\left\{\int_0^T\tau^{\theta p'}\|\partial_\tau v(\tau
)\|_{L^p(R)}^{p'}d\tau \right\}^{\frac{1}{p'}}\le
C(\|u_0\|_{L^p}+\int_0^T\|\partial_\tau w(\tau
)\|_{L^p(R)}d\tau)^3.
\end{equation}
By H${\rm \ddot{o}}$lder's inequality, we obtain
\begin{eqnarray}
&&\int_0^T\|\partial_\tau w(\tau )\|_{L^p}d\tau\le \left\{
\int_0^T \tau^{-\theta
p}\right\}^{\frac{1}{p}}\left\{\int_0^T\tau^{\theta
p'}\|\partial_\tau w(\tau
)\|_{L^p(R)}^{p'}d\tau\right\}^{\frac{1}{p'}}\\\nonumber &&
=CT^{\frac{1}{p'}}\left\{\int_0^T\tau^{\theta p'}\|\partial_\tau
w(\tau )\|_{L^p(R)}^{p'}d\tau\right\}^{\frac{1}{p'}}\\\nonumber &&
\le CC_1T^{\frac{1}{p'}}\|u_0\|_{L^p(R)}^3.
\end{eqnarray}
It then follows that
\begin{eqnarray}
&& \left\{\int_0^T\tau^{\theta p'}\|\partial_\tau v(\tau
)\|_{L^p(R)}^{p'}d\tau\right\}^{\frac{1}{p'}}\le
C(\|u_0\|_{L^p(R)}+C_1T^{\frac{1}{p'}}\|u_0\|_{L^p(R)}^3)^3\\\nonumber
&& \le C_1\|u_0\|_{L^p(R)}^3
\end{eqnarray}
provided that T is sufficiently small. By a similar argument, we
can show that M is a contraction. Moreover, it is not difficulty
to prove \eqref{a29}  and \eqref{a30}.

\section{Acknowledgement}
The author was supported by the National Natural Science
Foundation of China under grant 10225102.



\begin{thebibliography}{00}
\bibitem{Bergh} J. Bergh and J. L${\rm \ddot{o}}$fstr${\rm \ddot{o}}$m,
Interpolation spaces. An introduction. Grundlehren der
Mathematischen Wissenschaften, No. 223. Springer-Verlag,
Berlin-New York, 1976.

\bibitem{Bourgain} J. Bourgain, Fourier transform restriction
phenomena for certain lattice subsets and applications to
nonlinear evolution equations. I. Schr${\rm \ddot{o}}$dinger
equations. \textit{Geom. Funct. Anal.} \textbf{3} (1993), no. 2,
107--156.

\bibitem{Cazenave et al}T.Cazenave, L.Vega and M.C.Vilela, A note
on the nonlinear Schr${\rm \ddot{o}}$dinger equation in weak $L^p$
spaces. \textit{Communications in contemporary
Mathematics},\textbf{3}(2001), no.1,153-162.

\bibitem{CazenaveWeissler} T. Cazenave and F. B. Weissler,
The Cauchy problem for the critical nonlinear Schr${\rm
\ddot{o}}$dinger equation. \textit{Non. Anal. TMA}, \textbf{14}
(1990), 807-836.

\bibitem{Christ et al}M.Christ,J.Colliander and T.Tao, A priori
bounds and weak solutions for the nonlinear Schr${\rm
\ddot{o}}$dinger equation in Sobolev spaces of negative order.
\textit{preprint}.

\bibitem{Gruenrock} A.Grunr${\rm \ddot{o}}$ck, Bi-and trilinear
Schr${\rm \ddot{o}}$dinger estimate in one space dimension with
applications to cubic NLS and DNLS \textit{Int.Math.Res.Not.}
(2005), no.41, 2524-2558.


\bibitem{Klainerman1} S. Klainerman, Remarks on the global
Sobolev inequalities in the Minkowski space $R\sp {n+1}$.
\textit{Comm. Pure Appl. Math.} \textbf{40} (1987), no. 1,
111--117.

\bibitem{KochTataru}H.Koch and D.Tataru, A-priori bounds  for the
1-D cubic NLS in negative Sobolev spaces. \textit{preprint }.

\bibitem{McKeanShatah} H. P. McKean and J. Shatah, The nonlinear
Schr${\rm \ddot{o}}$dinger equation and the nonlinear heat
equation reduction to linear form. \textit{Comm. Pure Appl. Math.}
\textbf{44} (1991), no. 8-9, 1067--1080.


\bibitem{Tsutsumi} Y. Tsutsumi, $L^2$ solutions for nonlinear
Schr${\rm \ddot{o}}$dinger equations and nonlinear groups.
\textit{Funk. Ekva.} \textbf{30} (1987), 115-125.

\bibitem{VargasVega}A.Vargas and L.Vega, Global wellposedness for
1D non-linear Schr${\rm \ddot{o}}$dinger equation for data with an
infinite $L^2$ norm.\textit{J.Math.Pures Appl.} \textbf{80}
(2001), 1029-1044.


\end{thebibliography}
\end{document}